\documentclass[11pt]{article}
\usepackage{amsmath}
\usepackage{theorem}
{\theoremstyle{break}
\theorembodyfont{\it}
\theoremheaderfont{\normalfont \bfseries}
\newtheorem{prop}{Proposition}[section]
\newtheorem{thm}{Theorem}[section]
\newtheorem{cor}{Corollary}[section]
\newtheorem{lemma}{Lemma}[section]
}
{\theoremstyle{break}
\theorembodyfont{\normalfont}

}
{\theoremstyle{plain}
\newtheorem{defn}{Definition}[section]
\theoremheaderfont{\normalfont\bfseries}
}
\begin{document}

\newcommand{\ls}[1]
   {\dimen0=\fontdimen6\the\font \lineskip=#1\dimen0
\advance\lineskip.5\fontdimen5\the\font \advance\lineskip-\dimen0
\lineskiplimit=.9\lineskip \baselineskip=\lineskip
\advance\baselineskip\dimen0 \normallineskip\lineskip
\normallineskiplimit\lineskiplimit \normalbaselineskip\baselineskip
\ignorespaces }


\newcommand{\be}{\begin{equation}}
\newcommand{\ee}{\end{equation}}
\newcommand{\bea}{\begin{eqnarray}}
\newcommand{\eea}{\end{eqnarray}}
\newcommand{\beaa}{\begin{eqnarray*}}
\newcommand{\eeaa}{\end{eqnarray*}}
\newcommand{\beno}{\begin{eqalignno}}
\newcommand{\eeno}{\end{eqalignno}}
\newcommand{\betwo}{\begin{eqaligntwo}}
\newcommand{\eetwo}{\end{eqaligntwo}}

\newcommand{\ben}{\begin{enumerate}}
\newcommand{\een}{\end{enumerate}}
\newcommand{\bi}{\begin{itemize}}
\newcommand{\ei}{\end{itemize}}

\newcommand{\lip}{\langle}
\newcommand{\rip}{\rangle}
\newcommand{\uu}{\underline}
\newcommand{\oo}{\Omega}
\newcommand{\La}{\Lambda}
\newcommand{\la}{\lambda}
\newcommand{\eps}{\epsilon}
\newcommand{\veps}{\varepsilon}
\newcommand{\om}{w}

\newcommand{\bige}{\mbox{\Large\it e}}
\newcommand{\integers}{Z\!\!\!Z}
\newcommand{\rationals}{{\rm I\!Q}}
\newcommand{\ID}{{\rm I\!D}}
\newcommand{\IN}{{\rm I\!N}}
\newcommand{\reals}{{\rm I\!R}}
\newcommand{\realsd}{\reals^d}
\newcommand{\degree}{{\scriptscriptstyle \circ }}
\newcommand{\dfn}{\stackrel{\triangle}{=}}
\def\complex{\mathop{\raise .45ex\hbox{${\bf\scriptstyle{|}}$}
     \kern -0.40em {\rm \textstyle{C}}}\nolimits}
\newcommand{\NN}{{\rm I\!N}}
\newcommand{\R}{{\rm I\!R}}
\newcommand{\RAISE}{{\:\raisebox{.8ex}{$\scriptstyle{>}$}\raisebox{-.3ex}
           {$\scriptstyle{\!\!\!\!\!\sim}\:$}}}
\newcommand{\Gi}{{\:\raisebox{.0ex}{$\textstyle{\cal X}$}\raisebox{.2ex}
           {$\textstyle{\!\!\!\!\! -}\:$}}}
\newcommand{\gi}{{\:\raisebox{.0ex}{$\scriptscriptstyle{\cal
        X}$}\raisebox{.1ex} 
           {$\scriptstyle{\!\!\!\!-}\:$}}}

\def\le{\leq}
\def\ge{\geq}
\def\lt{<}
\def\gt{>}

\newcommand{\loc}{{\rm loc}}
\newcommand{\Dom}{{\rm Dom}}
\newcommand{\calM}{{\cal M}}
\newcommand{\calA}{{\cal A}}
\newcommand{\calB}{{\cal B}}
\newcommand{\calC}{{\cal C}}
\newcommand{\calD}{{\cal D}}
\newcommand{\calF}{{\cal F}}
\newcommand{\calL}{{\cal L}}
\newcommand{\calP}{{\cal P}}
\newcommand{\calX}{{\cal X}}
\newcommand{\calH}{{\cal H}}

\newcommand{\DD}{{\rm I\!D}}
\newcommand{\dett}{{\mbox{\rm det}}_2}
\newcommand{\trace}{{\mbox{\rm trace}}}
\newcommand{\Prob}{{\rm Prob\,}}
\newcommand{\sinc}{{\rm sinc\,}}
\newcommand{\ctg}{{\rm ctg\,}}
\newcommand{\ifff}{\mbox{\ if and only if\ }}
\newcommand{\proof}{\noindent {\bf Proof:\ }}
\newcommand{\remark}{\noindent {\bf Remark:\ }}
\newcommand{\remarks}{\noindent {\bf Remarks:\ }}
\newcommand{\note}{\noindent {\bf Note:\ }}

\newcommand{\boldx}{{\bf x}}
\newcommand{\boldX}{{\bf X}}
\newcommand{\boldy}{{\bf y}}
\newcommand{\uux}{\uu{x}}
\newcommand{\uuY}{\uu{Y}}

\newcommand{\al}{\alpha}

\newcommand{\limn}{\lim_{n \rightarrow \infty}}
\newcommand{\limN}{\lim_{N \rightarrow \infty}}
\newcommand{\limr}{\lim_{r \rightarrow \infty}}
\newcommand{\limd}{\lim_{\delta \rightarrow \infty}}
\newcommand{\limM}{\lim_{M \rightarrow \infty}}
\newcommand{\limsupn}{\limsup_{n \rightarrow \infty}}

\newcommand{\imii}{\int_{-\infty}^{\infty}}
\newcommand{\imix}{\int_{-\infty}^x}
\newcommand{\ioi}{\int_o^\infty}

\newcommand{\ARROW}[1]
  {\begin{array}[t]{c}  \longrightarrow \\[-0.2cm] \textstyle{#1} \end{array} }

\newcommand{\AR}
 {\begin{array}[t]{c}
  \longrightarrow \\[-0.3cm]
  \scriptstyle {n\rightarrow \infty}
  \end{array}}

\newcommand{\pile}[2]
  {\left( \begin{array}{c}  {#1}\\[-0.2cm] {#2} \end{array} \right) }

\newcommand{\floor}[1]{\left\lfloor #1 \right\rfloor}

\newcommand{\mmbox}[1]{\mbox{\scriptsize{#1}}}

\newcommand{\ffrac}[2]
  {\left( \frac{#1}{#2} \right)}

\newcommand{\one}{\frac{1}{n}\:}
\newcommand{\half}{\frac{1}{2}\:}

\def\squarebox#1{\hbox to #1{\hfill\vbox to #1{\vfill}}}
\newcommand{\qed}{\hspace*{\fill}
            \vbox{\hrule\hbox{\vrule\squarebox{.667em}\vrule}\hrule}\smallskip}

\thispagestyle{empty}

\vspace*{.5cm}
\begin{center}

{\LARGE\bf The Notion of Convexity and Concavity on Wiener Space}\\

\vspace*{.5cm}
D. Feyel and A. S.  \"{U}st\"unel
\end{center}

\vspace*{0.5cm}
\noindent
{\bf Abstract}--
{\footnotesize {We define, in the frame of an abstract Wiener space,  the  notions  of
convexity and of  concavity for the equivalence classes of random
variables. As application we show that some important inequalities of
the finite dimensional case have their natural counterparts in this setting.}}
\vspace*{0.5cm}
\noindent
\section{Introduction}
On an infinite dimensional vector space $W$  the notion of convex or concave
function is well-known. Assume now that this space is equipped with a
probability measure. Suppose that there are two measurable functions  on 
this vector space, say $F$ and $G$ such that $F=G$ almost surely. If
$F$ is a  convex function, then from the probabilistic point of view,
we would like to say  that $G$ is also convex. However this is false;  since
in general the underlying probability measure is not (quasi) 
invariant under the translations by the elements of the vector
space. If $W$ contains a dense  subspace $H$ such that $w\to w+h$
($h\in H$) induces a measure which is equivalent to the initial
measure or absolutely continuous with respect to it, then   we can
define a  notion of ``$H$--convexity'' or ``$H$--concavity  in the
direction  of $H$. Of course  these  properties  are  inherited by the
corresponding equivalence classes, hence they are particularly useful
for the probabilistic calculations. 

The notion of $H$-convexity  has
been used  in \cite{USZ96} to study the 
absolute continuity of the image of the Wiener measure under the 
monotone shifts. In this paper we study further  properties of  
such functions  and some additional ones  in the frame of an abstract
Wiener space, namely  $H$-convex, $H$-concave, log $H$-concave and log 
$H$-convex Wiener functions,
where $H$ denotes the associated  Cameron-Martin space. In particular
we extend some finite dimensional results of  \cite{PRE71} and
\cite{BLI}  to this setting and prove that some finite dimensional
convexity-concavity inequalities have their counterparts in infinite
dimensions. 
\section{Preliminaries}

In the sequel $(W,H,\mu)$ denotes an abstract Wiener space, i.e., $H$
is a separable Hilbert space, called the Cameron-Martin space. 
It is identified with its continuous dual.  $W$ 
is a Banach or a Fr\'echet space into which $H$ is injected
continuously and densely. $\mu$ is the standard  cylindrical  Gaussian
measure on $H$ which is concentrated in $W$ as a Radon probability measure.
In the  classical case  we have either $ W= C_0([0,1]) $ or
$W=C_0(\reals_+)$ where 
$$ 
H=\left\{h : [0,1]\rightarrow \reals : h(t) = \int_0^t \, \dot{h} (s) \, ds\,, 
|h|_H = \|\dot{h}\|_{L^2([0,1])}\right\}\,
$$
or 
$$ 
H=\left\{h :\reals_+\rightarrow \reals : h(t) = \int_0^t \, \dot{h}
  (s)\, ds\,, |h|_H = \|\dot{h}\|_{L^2(\reals_+)}\right\}\, 
$$
respectively.

Let $X$ be a separable Hilbert space and $a$ be an $X$-valued (smooth)
polynomial on $W$:
$$
a(w) = \sum_{i=1}^m\, \eta_i (\langle h_1 , w\rangle \, , \ldots , \,
\langle h_n , w\rangle ) x_i \; ,
$$
with $ x_i \in X , \, h_i \in W^* $ and
$ \eta_i \in C_b^{\infty} (\reals^n ) $.
The Gross-Sobolev derivative of $a$ is defined as
$$
\nabla a (w) = \sum_{i=1}^m\, \sum_{j=1}^n \,
\partial_j \eta_i (\lip h_1 , w \rip \, , \ldots , \, \lip h_n , w\rip ) x_i
\otimes {\tilde {h}}_j \; ,
$$
where $\tilde{h}$ denotes the image of $h\in W^*$ in $H$ under the
canonical injection $W^*\hookrightarrow H$ (in the sequel we shall
omit this notational detail and write $h$ instead of $\tilde{h}$ when
there is no ambiguity).
The derivatives of higher orders  $\nabla^k a(w) $ are  defined recursively.
Thanks to the Cameron-Martin theorem, all these operators are closable on all
the $ L^p $--spaces and the Sobolev spaces
$ \DD_{p,k} (X) , \, p\gt 1 , \; k\in\NN $
can be defined as the completion of $X$-valued smooth polynomials with
respect to the norm:
$$
\parallel a\parallel_{p,k} =\sum_{i=0}^k \, \parallel\nabla^i  a\parallel_{L^p
(\mu , X \otimes H^{\otimes i} ) } \; .
$$
 From the Meyer inequalities (cf., for instance \cite{UST95}), it is known that
 the $ (p,k)$-norm, defined above, is equivalent to the following norm
$$
\parallel (I+L)^{k/2} a\parallel_{L^p (\mu, X)}
$$
where $L$ is the Ornstein-Uhlenbeck operator on $W$ (cf. \cite{UST95})
and we denote these two norms with the same notation.
Since $L$ is a positive, self adjoint  operator, we can also define the norms,
via spectral theorem, for  $k\in \reals $. It is easy to see that the spaces
with negative differentiability index   describe the dual spaces of the
positively indexed Sobolev spaces.
We denote by $ \DD \, (X) $ the intersection of the Sobolev spaces
$ \{ \DD_{p,k} (X) ; \; p\gt 1 , \; k\in \NN\} $, equipped with the
intersection (i.e., projective limit) topology.
The continuous dual of $ \DD (X) $ is denoted by $ \DD'(X) $ and in case 
$X=\reals $ we write simply $ \DD_{p,k} , \, \DD , \, \DD' $ for
$ \DD_{p,k} (\reals) , \, \DD (\reals) , \, \DD' (\reals) $ respectively.
Consequently, for any $ p\gt 1 , \, k\in \reals , \; \nabla : \DD_{p,k} (X)
\mapsto \DD_{p,k-1} (X\otimes  H) $ continuously, where
$ X\otimes H$ denotes the completed Hilbert-Schmidt tensor product of
$X$ and $H$.
Therefore $\delta = \nabla^* $ is a continuous operator from $ \DD_{p,k}
(X\otimes H) $ into $ \DD_{p, k-1} (X) $ for any $ p\gt 1 , \; k\in \reals$.
We call $\delta$ the divergence operator on $W$. Let us remark that
from these  properties, $\delta$ and $\nabla$ extend continuously as
operators from $\DD'(X\otimes H)$ to $\DD'(X)$ and from $\DD'(X)$ to
$\DD'(X\otimes H)$ respectively.
Let us recall that, in the case of classical Wiener space, $\delta$ coincides
with the It\^{o} stochastic integral on the adapted processes.
We recall that, if $F$ is in $\DD_{p,1} (H) $ for some $ p \gt 1 $, then almost
surely, $\nabla F $ is an Hilbert-Schmidt operator on $H$, and if $F$ is an
$H$-valued polynomial, then $\delta F $ can be written as
$$
\delta F =\sum_{i=1}^{\infty} \, \Bigl[(F, e_i)_H\delta e_i - \Bigl(\nabla (F, e_i)_H
, e_i\Bigr)_H \Bigr] \; ,
$$
where $ (e_i,\,i\in \NN ) $ is any complete orthonormal basis in $H$.


In the sequel we shall use the notion of second quantization of 
bounded operators on $H$; although this is a well-known subject, we
give a brief outline below for the reader's convenience
(cf. \cite{BAD}, \cite{P-F}, \cite{SIM74}). 
Assume that $A:H\to H$ is a bounded, linear operator, then it has a
unique, $\mu$-measurable (i.e., measurable with respect to the
$\mu$-completion of $\calB(W)$) extension, denoted by $\tilde{A}$,  as
a linear map on $W$ (cf.\cite{BAD,P-F}). Assume in particular  
that $\|A\|\leq 1$ and define $S=(I_H-A^*A)^{1/2}$,
$T=(I_H-AA^*)^{1/2}$ and $U:H\times H\to H\times H$ as
$U(h,k)=(Ah+Tk,-Sh+A^*k)$. $U$ is then  a  unitary operator on $H\times H$,
hence its $\mu\times\mu$-measurable linear extension to $W\times W$
preserves the Wiener measure $\mu\times\mu$ (this is called the
rotation associated 
to $U$, cf. \cite{USZ99}, Chapter VIII). Using this observation, one can define
the second quantization of $A$ via the generalized Mehler formula as  
$$
\Gamma(A)f(w)=\int_Wf(\tilde{A^*}w+\tilde{S}y)\mu(dy)\,, 
$$
which happens to be a Markovian contraction on $L^p(\mu)$ for any
$p\geq 1$. $\Gamma (A)$ can be calculated explicitly for the Wick
exponentials as 
$$
\Gamma(A)\exp\left\{\delta h-1/2|h|_H^2\right\}=\exp\left\{\delta A
  h-1/2|Ah|_H^2\right\}\,\,(h\in H)\,. 
$$
This identity implies that $\Gamma (AB)=\Gamma(A)\Gamma(B)$ and that
for any sequence $(A_n,n\in \NN)$ of operators whose norms are bounded 
by one, $\Gamma(A_n)$ converges strongly to $\Gamma(A)$ if
$\lim_nA_n=A$ in the strong operator topology.
 A particular case of interest is  when we take $A=e^{-t}I_H$, then
 $\Gamma(e^{-t}I_H)$ equals to  the Ornstein-Uhlenbeck semigroup
 $P_t$. Also  if 
 $\pi$ is the  orthogonal projection of $H$ onto a closed vector subspace
 $K$, then $\Gamma(\pi)$ is the conditional expectation with respect
 to the sigma field generated  by $\{\delta k,\,k\in K\}$.



\section{$H$-convexity and its properties}
Let us give  the notion of $H$-convexity on the  Wiener space $W$:
\begin{defn}
Let  $F:W\to \reals\cup\{\infty\}$ be a measurable function. 
It is called $H$-convex if for any 
$h,k\in H$, $\al\in [0,1]$
\begin{equation}
\label{con-def}
F(w+\al h+(1-\al)k)\leq \al F(w+h)+(1-\al)F(w+k)
\end{equation}
almost surely.
\end{defn}
\remarks
\begin{itemize}
\item This definition is more general than the one given in
  \cite{USZ96,USZ99} since $F$ may be   infinite on a set of  positive 
  measure.
\item Note that the negligeable  set on which the relation
  (\ref{con-def}) fails may depend on the choice of $h,k$ and of $\al$.
\item If $G:W\to\reals$ is a measurable  convex function, then it is
  necessarily $H$-convex.
\item To conclude the $H$-convexity, it suffices to verify the
  relation (\ref{con-def}) for $k=-h$ and $\alpha=1/2$.
\end{itemize}

\noindent
The following properties of $H$-convex Wiener functionals have been
proved in  \cite{USZ96,USZ98,USZ99}:
\begin{thm}
\label{char-thm}
\begin{enumerate}
\item If $(F_n,n \in \NN)$ is a sequence  of $H$-convex 
functionals converging in probability, then the limit is also
  $H$-convex.
\item If $F\in L^p(\mu)$ ($p>1$) is $H$-convex if and only if
  $\nabla^2 F$ is positive and symmetric  Hilbert-Schmidt operator
  valued distribution on $W$.
\item If $F\in L^1(\mu)$ is $H$-convex, then $P_tF$ is also $H$-convex 
  for any $t\geq 0$, where $P_t$ is the Ornstein-Uhlenbeck semi-group
  on $W$.
\end{enumerate}
\end{thm}
\noindent
The following result is immediate from Theorem \ref{char-thm} :
\begin{cor}
\label{char-cor}
$F\in \cup_{p>1}L^p(\mu)$ is $H$-convex if and only if 
$$
E\left[\varphi\,\left(\nabla^2F(w),h\otimes h\right)_2\right]\geq 0
$$
for any $h\in H$ and $\varphi\in \DD_+$, where $(\cdot\,,\,\cdot)_2$
denotes the scalar 
product for the Hilbert-Schmidt operators on $H$ .
\end{cor}
We have also
\begin{cor}
If $F\in L^p(\mu)$, $p>1$, is $H$-convex and if $E[\nabla^2F]=0$, then 
$F$ is of the form
$$
F=E[F]+\delta\left(E[\nabla F]\right)\,.
$$
\end{cor}
\proof
Let $(P_t,t\geq 0)$ denote the Ornstein-Uhlenbeck semigroup, $P_tF$ is
again $H$-convex and Sobolev differentiable.  Moreover $\nabla^2
P_tF=e^{-2t}P_t\nabla^2F$. Hence $E[\nabla^2P_tF]=0$, and  the
positivity of $\nabla^2P_tF$ implies that $\nabla^2 P_tF=0$ almost
surely, hence $\nabla^2 F=0$. This implies that $F$ is in the first
two Wiener chaos.
\qed

\remark It may be worth-while  to note that the random variable which
represents  the share price of the Black and
Scholes model in financial mathematics (cf.\cite{L-L})  is $H$-convex.
\newline

We shall need also the  concept of $\calC$-convex functionals: 
\begin{defn}
\label{C-defn}
 Let $(e_i,i\in \NN)\subset W^*$ be any complete, orthonormal basis of
 $H$. For $w\in W$, define  $w_n=\sum_{i=1}^n\delta e_i(w)e_i$ and
 $w_n^\perp=w-w_n$, then a  Wiener functional $f:W\to \R$ is called
 $\calC$-convex if, for any such basis $(e_i,i\in \NN)$, for almost
 all $w_n^\perp$,  the partial map 
$$
w_n\to f(w_n^\perp+w_n)
$$
has a modification which is convex on the space 
${\rm span}\{e_1,\ldots,e_n\}\simeq  \R^n$. 
\end{defn}
\remark It follows from Corollary \ref{char-cor}  that, if $f$ is
$H$-convex and   in
some $L^p(\mu)$ $(p>1)$, then it is $\calC$-convex. We shall prove
that this is also true without any integrability hypothesis.

We begin with the following lemma whose proof is obvious:

\begin{lemma}
\label{L1}
If $f$ is $\calC$-convex then it is $H$-convex.
\end{lemma}
\noindent
In order to  prove the validity of the  converse of Lemma
\ref{L1}  we need some technical results from  the
harmonic analysis on finite dimensional Euclidean spaces that we shall 
state as separate lemmas:
\begin{lemma}
\label{L2}
Let $B\in {\calB}(\reals^n)$ be a set of positive Lebesgue
measure. Then $B+B$ contains a non-empty open set.
\end{lemma}
\proof
Let $\phi(x)=1_B\star 1_B(x)$, where ``$\star$''  denotes the
convolution of functions  with respect to the Lebesgue measure. Then
$\phi$ is a non-negative,  continuous function, hence the set $O=\{x\in
\reals^n:\,\phi(x)>0\}$ is an open set. Since $B$ has positive
measure, $\phi$ can not be identically zero, hence $O$ is
non-empty. Besides, if $x\in O$, then the set of $y\in \reals^n$ such
that $y\in B$ and $x-y\in B$ has positive Lebesgue measure, otherwise
$\phi(x)$ would have been null. Consequently $O\subset B+B$.
\qed

The following lemma gives a more precise statement than Lemma \ref{L2}:  
\begin{lemma}
\label{L3}
Let $B\in {\calB}(\reals^n)$ be a set of positive Lebesgue
measure and  assume  that  $A\subset \reals^n\times \reals^n$ with 
$B\times B=A$ almost surely with respect to the Lebesgue measure of
$\reals^n\times \reals^n$. Then the set $\{x+y:\,(x,y)\in A\}$
contains almost surely an open subset of $\reals^n$.
\end{lemma}
\proof
It follows from an obvious change of variables that 
$$
1_A(y,x-y)=1_B(y)1_B(x-y)
$$ 
almost surely, hence 
$$
\int_{\reals^n}1_A(y,x-y)dy=\phi(x)
$$
almost surely, where $\phi(x)=1_B\star 1_B(x)$. Consequently, for
almost all $x\in \reals^n$ such that $\phi(x)>0$, one has $(y,x-y)\in
A$, this means that 
$$
\{x\in \reals^n:\,\phi(x)>0\}\subset
\{u+v:\,(u,v)\in A\}
$$ 
almost surely.
\qed 

The following lemma is particularly important for the sequel:
\begin{lemma}
\label{L4}
Let $f:\reals^n\to \reals_+\cup\{\infty\}$ be a Borel function which
is finite on a set of positive Lebesgue measure. Assume that, for any
$u\in \reals^n$, 
\begin{equation}
\label{C1}
f(x)\leq \frac{1}{2}[f(x+u)+f(x-u)]
\end{equation}
$dx$-almost surely (the negligeable set on which the inequality
(\ref{C1}) fails may depend on $u$). Then there exists a non-empty,
open convex subset $U$ of $\reals^n$ such that $f$ is locally essentially 
bounded on $U$. Moreover let $D$ be the set consisting of  
 $x\in \reals^n$ such  that  any neighbourhood
 of $x\in D$ contains a Borel  set of positive Lebesgue measure on
 which $f$ is finite, then $D\subset \overline{U}$, in particular
 $f=\infty$ almost surely  on the complement of $\overline{U}$. 
\end{lemma}
\proof
From the theorem of Fubini, the inequality (\ref{C1}) implies that 
\begin{equation}
\label{C2}
2f\left(\frac{x+y}{2}\right)\leq f(x)+f(y)
\end{equation}
$dx\times dy$-almost surely. Let $B\in {\calB}(\reals^n)$ be a set of
positive Lebesgue measure on which $f$ is  bounded by some constant
$M>0$. Then from Lemma 
\ref{L2}, $B+B$ contains an open set $O$. Let $A$ be the set
consisting of the elements of $B\times B$ for which the inequality
(\ref{C2}) holds. Then $A=B\times B$ almost surely, hence from Lemma
\ref{L3}, the set $\Gamma=\{x+y:\,(x,y)\in A\}$ contains almost surely 
the open set $O$. Hence for almost all $z\in \frac{1}{2}O$, $2z$
belongs to the set $\Gamma$, consequently $z=\frac{1}{2}(x+y)$, with
$(x,y)\in A$. This implies, from (\ref{C2}),  that $f(z)\leq
M$. Consequently $f$ is essentially bounded on the open set
$\frac{1}{2}\Gamma$. 

Let now $U$ be set of points  which have neighbourhoods on which $f$ is 
essentially bounded. Clearly $U$ is open and non-empty by what we have 
shown above. Let $S$ and $T$ be two balls of radius $\rho$, on which
$f$ is bounded by some $M>0$. Assume that they are centered at the
points $a$ and $b$ respectively. Let $u=\frac{1}{2}(b-a)$, then for
almost all $x\in \frac{1}{2}(S+T)$, $x+u\in T$ and $x-u\in S$, hence,
from the inequality (\ref{C1})
$f(x)\leq M$, which shows that $f$ is essentially bounded on the set
$\frac{1}{2}(S+T)$ and this proves the convexity of $U$. 

To prove the
last  claim, let $x$ be any element of $D$ and let $V$ be any
neighbourhood of $x$; without loss of generality, we may assume that
$V$ is convex. Then there exists a Borel set $B\subset V$ of positive
measure on which $f$ is bounded, hence from the first part of the
proof, there exists an open  neighbourhood $O\subset B+B$ such that $f$ is
essentially bounded on $\frac{1}{2}O\subset \frac{1}{2}(V+V)\subset
V$, hence $\frac{1}{2}O\subset U$. Consequently $V\cap U\neq
\emptyset$, and this implies that $x$ is in the closure of $U$,
i.e. $D\subset \overline{U}$. The fact that $f=\infty$ almost surely
on the complement of $\overline{U}$ is obvious from the definition of $D$.
\qed

\begin{thm}
\label{con-char}
Let $g:\reals^n\to \reals\cup\{\infty\}$ be a measurable mapping such
that, for almost all $u\in \reals^n$,  
\begin{equation}
\label{CC}
g(u+\alpha x+\beta y)\leq \alpha g(u+x)+\beta g(u+y)
\end{equation}
for any $\alpha,\beta \in [0,1]$ with $\alpha+\beta =1$ and for any
$x,y\in \R^n$, where the
negligeable set on which the relation (\ref{CC}) fails may depend on
the choice of $x,y$ and of $\alpha$ . Then $g$ has
a modification $g'$ which is a convex function.
\end{thm}
\proof
Assume first that $g$ is positive, then with the notations of Lemma
\ref{L4}, define $g'=g$ on the open, convex set $U$ and as $g'=\infty$ 
on $U^c$. From the relation (\ref{CC}), $g'$ is a distribution on $U$
whose second derivative is positive, hence it is convex on $U$, hence
it is convex on the whole space $\reals^n$. Moreover we have $\{g'\neq 
g\}\subset \partial U$ and $\partial U$ has zero Lebesgue measure,
consequently $g=g'$ almost surely. For general $g$, define
$f_\eps=e^{\eps g}$ ($\eps>0$), then, from what is proven above,
  $f_\eps$ has a modification $f'_\eps$ which is convex (with the same 
  fixed open and convex set $U$), hence $\lim\sup_{\eps\to
    0}\frac{f'_\eps-1}{\eps}=g'$ is also convex and $g=g'$ almost
  surely.
\qed

\begin{thm}
\label{H->C}
A Wiener functional $F:W\to \reals\cup\{\infty\}$ is $H$-convex if and 
only if it is $\calC$-convex.
\end{thm}
\proof
We have already proven the sufficiency. To prove the necessity, with
the notations of Definition \ref{C-defn}, $H$-convexity implies that 
$h\to F(w_n^\perp+w_n+h)$ satisfies the hypothesis of Theorem
\ref{con-char} when $h$ runs in any $n$-dimensional Euclidean subspace 
of $H$, hence the partial mapping $w_n\to F(w_n^\perp+w_n)$ has a
modification which is convex on the vector space spanned by
$\{e_1,\ldots,e_n\}$. 
\qed

\section{Log $H$-concave and $\calC$-$\log$ concave  Wiener functionals}
 
\begin{defn}
Let $F$ be a  measurable mapping  from $W$ into $\reals_+$ with
$\mu\{F>0\}>0$.  
\begin{enumerate}
\item $F$ is  called log $H$-concave, if for any $h,\,k\in H$,
  $\alpha\in [0,1]$, one has  
\begin{equation}
\label{log-cv-def}
F\bigl(w+\alpha h+(1-\alpha)k\bigr)\geq F(w+ h)^\alpha\,  F(w+k)^{1-\alpha}
\end{equation}
almost surely, where the negligeable set on which the relation
(\ref{log-cv-def}) 
fails may depend on $h,\,k$ and on $\alpha$.
\item We shall say that $F$ is $\calC$-log concave, if for any
  complete, orthonormal basis  $(e_i,i\in \NN)\subset W^*$ of $H$, the 
  partial map $w_n\to F(w_n^\perp+w_n)$ is log-concave (cf. Definition 
  \ref{C-defn} for the notation),  up to a
  modification,  on ${\rm span}\{e_1,\ldots,e_n\}\simeq \R^n$.
\end{enumerate}
\end{defn}

Let us  remark immediately that if $F=G$ almost surely then $G$ is also log
$H$-concave. Moreover, any limit in probability of log $H$-concave random
variables is again log $H$-concave. We shall prove below some less immediate
properties. Let us begin with the following observation which is a
direct consequence of Theorem \ref{H->C}:

\remark
\label{H-concave=C-concave}
$F$ is log $H$-concave if and only if $-\log F$ is $H$-convex 
 (which may be infinity with a positive probability), hence if and
 only if $F$ is $\calC$-log concave.

\begin{thm}
\label{prekopa1-the}
Suppose that $(W_i,H_i,\mu_i)$, $i=1,2$, are two abstract Wiener
spaces. Consider $(W_1\times W_2, H_1\times H_1,\mu_1\times \mu_2)$ as an
abstract Wiener space. Assume that $F:W_1\times W_2\to \reals_+$ is log
$H_1\times H_2$-concave. Then the map 
$$
w_2\to \int_{W_1}F(w_1,w_2)\,d\mu_1(w_1)
$$
is log $H_2$-concave.
\end{thm}
\proof
If $F$ is log $H\times H$-concave, so is also $F\wedge c$ ($c\in
\reals_+$), hence  we may suppose without
loss of generality that $F$ is bounded.
Let $(e_i,i\in \NN)$ be a complete, orthonormal basis in $H_2$. It suffices to
prove that 
$$
E_1[F](w_2+\alpha h+\beta l)\geq\left(E_1[F](w_2+h)\right)^\alpha
  \left(E_1[F](w_2+l)\right)^\beta
$$
almost surely, for any $h,\,l\in {\mbox{span}}\{e_1,\ldots,e_k\}$,
$\alpha,\beta\in [0,1]$ with $\alpha+\beta=1$, where $E_1$ denotes the
expectation with respect to $\mu_1$. Let $(P_n,n\in \NN)$ be a sequence of
orthogonal projections of finite rank  on $H_1$ increasing to the
identity map of it. Denote by 
$\mu_1^n$ the image of $\mu_1$ under the map $w_1\to {\tilde{P}}_n w_1$ and by
$\mu_1^{n\perp}$ the image of $\mu_1$ under $w_1\to w_1-{\tilde{P}}_n w_1$. We
have, from the martingale convergence theorem,  
$$
\int_{W_1}F(w_1,w_2)\,d\mu_1(w_1)=\lim_n\int
F(w_1^{n\perp}+w_1^n,w_2)\,d\mu_1^n(w^n_1)
$$
almost surely. Let $(Q_n,n\in \NN)$ be a sequence of orthogonal
projections of finite rank  on
$H_2$ increasing to the identity, corresponding to the basis $(e_n,n\in
\NN)$. Let $w_2^k={\tilde{Q}}_kw_2$ and $w_2^{k\perp}=w_2-w_2^k$. Write 
\beaa
F(w_1,w_2)&=&F(w_1^{n\perp}+w_1^n,w_2^k+w_2^{k\perp})\\
       &=&F_{w_1^{n\perp},w_2^{k\perp}}(w_1^n,w_2^k)\,.
\eeaa
From the hypothesis
$$
(w_1^n,w_2^k)\to F_{w_1^{n\perp},w_2^{k\perp}}(w_1^n,w_2^k)
$$
has a  log concave modification  on the $(n+k)$-dimensional Euclidean
space. From the theorem of Pr\'ekopa (cf.\cite{PRE71}), it follows that 
$$
w_2^k\to \int F_{w_1^{n\perp},w_2^{k\perp}}(w_1^n,w_2^k)\,d\mu_1^n(w_1^n)
$$
is log concave on $\R^k$  for any $k\in \NN$ (upto a modification), hence 
$$
w_2\to \int F(w_1^{n\perp}+w_1^n,w_2)\, d\mu(w_1^n)
$$
is log $H_2$-concave for any $n\in \NN$, then the proof follows by passing to
the limit with respect to $n$.
\qed

\begin{thm}
\label{Gamma-thm}
Let $A:H\to H$ be a linear operator with $\|A\|\leq 1$, denote by
$\Gamma(A)$ its second  quantization as explained in the
preliminaries. If $F:W\to \reals_+$ is a log $H$-concave Wiener
functional, then  $\Gamma(A)F$  is also log $H$-concave.
\end{thm}
\proof
Replacing $F$ by $F\wedge c=\min(F,c),\,c>0$, we may suppose that $F$
is bounded. It is
easy  to see that the mapping 
$$
(w,y)\to F(\tilde{A^*}w+\tilde{S}y)
$$
is log $H\times H$-concave on $W\times W$. In fact, for any
$\alpha+\beta=1$, $h,k,u,v\in H$, one has 
\bea
\label{formul}
\lefteqn{F(\tilde{A^*}w+\tilde{S}y+\alpha(A^*h+Sk)+\beta(A^*u+Sv))}\\
&&\geq F(\tilde{A^*}w+\tilde{S}y+A^*h+Sk)^\alpha\,
F(\tilde{A^*}w+\tilde{S}y+A^*u+Sv)^\beta\,,\nonumber
\eea
$d\mu\times d\mu$-almost surely. Let us recall that, since the image
of $\mu\times \mu$ under the map $(w,y)\to \tilde{A^*}w+\tilde{S}y$ 
is $\mu$, the terms in the inequality (\ref{formul}) are defined without 
ambiguity. Hence 
$$
\Gamma(A)F(w)=\int_W F(\tilde{A^*}w+\tilde{S}y)\mu(dy)
$$
is log $H$-concave on $W$ from Theorem \ref{prekopa1-the}. 
\qed

\begin{cor}
\label{con-cor}
Let $F:W\to \reals_+$ be a log $H$-concave functional. Assume that $K$ 
is any closed vector subspace of $H$ and denote by $V(K)$ the sigma algebra
generated by $\{\delta k,\,k\in K\}$. Then the  conditional
expectation of $F$ with respect to $V(K)$, i.e., $E[F|V(K)]$  is again
log $H$-concave. 
\end{cor}
\proof
The proof follows from Theorem \ref{Gamma-thm} as soon as we remark
that $\Gamma(\pi_K)F=E[F|V(K)]$, where $\pi_K$ denotes the orthogonal
projection associated to $K$.
\qed

\begin{cor}
\label{semi-con}
Let $F$  be log $H$-concave. If $P_t$ denotes the
Ornstein-Uhlenbeck semigroup on $W$, then $w\to P_tF(w)$ is log $H$-concave.
\end{cor}
\proof
Since $P_t=\Gamma(e^{-t}I_H)$, the proof follows from Theorem
\ref{Gamma-thm}.
\qed

Here is an important application of these results:
\begin{thm}
\label{version-thm}
Assume that $F:W\to \reals\cup \{\infty\}$ is an $H$-convex Wiener
functional, then $F$ has a modification $F'$ which is a Borel
measurable  convex
function on $W$.  Any log $H$-concave 
functional $G$ has a modification $G'$ which is Borel measurable and
log-concave on $W$.  
\end{thm} 
\proof
Assume first that $F$ is positive, let $G=\exp-F$, then $G$ is a
positive, bounded $\calC$-log concave function. Define $G_n$ as  
$$
G_n=E[P_{1/n}G|V_n]\,,
$$
where $V_n$ is the sigma algebra generated by $\{\delta
e_1,\ldots,\delta e_n\}$, and $(e_i,i\in \NN)\subset W^*$ is a
complete orthonormal basis of $H$. Since
$P_{1/n}E[G|V_n]=E[P_{1/n}G|V_n]$, the positivity improving property
of the Ornstein-Uhlenbeck semigroup implies that $G_n$ is almost
surely strictly positive  (even quasi-surely). As we have attained the 
finite dimensional case, $G_n$ has a modification 
$G'_n$ which is continuous on $W$ and, from Corollary \ref{con-cor}
and Corollary \ref{semi-con}, it satisfies
\begin{equation}
\label{f1}
G'_n(w+ah+bk)\geq G'_n(w+h)^a G'_n(w+k)^b
\end{equation}
almost surely, for any $h,k\in H$ and $a+b=1$. The continuity of
$G'_n$ implies that  the relation (\ref{f1}) holds for any $h,k\in H$,
$w\in W$ and $a\in [0,1]$. Hence $G'_n$ is log-concave on $W$ and 
this implies that $-\log G'_n$ is convex on $W$. Define
$F'=\lim\sup_n(-\log G'_n)$, then $F'$ is convex and Borel measurable
on $W$ and $F=F'$ almost surely. 

For general $F$, define $f_\eps=e^{\eps F}$, then from above, there
exists a modification of $f_\eps$, say  $f'_\eps$ which is convex and
Borel measurable  on
$W$. To complete the proof it suffices to define $F'$ as  
$$
F'=\lim\sup_{\eps\to 0}\frac{f'_\eps-1}{\eps}\,.
$$
The rest is   now obvious.
\qed

Under the light of Theorem \ref{version-thm}, the following definition 
is natural:
\begin{defn}
A Wiener functional $F:W\to \reals\cup\{\infty\}$ will be called
almost surely convex if it has a modification $F'$ which is convex and 
Borel measurable on $W$. Similarly, a non-negative functional $G$ will
be called almost surely log-concave if it has a modification $G'$
which is log-concave on $W$.
\end{defn}

The following proposition summarizes the main results of this section:
\begin{thm}
\label{equiv-prop}
Assume that $F:W\to \reals\cup\{\infty\} $ is a Wiener functional such that
$$
\mu\{F<\infty\}>0\,.
$$
Then the following are equivalent:
\begin{enumerate}
\item $F$ is $H$-convex,
\item $F$ is $\calC$-convex,
\item $F$ is almost surely convex.
\end{enumerate}
Similarly, for $G:W\to \reals_+$, with $\mu\{G>0\}>0$, the following
properties are equivalent:
\begin{enumerate}
\item $G$ is log $H$-concave,
\item $G$ is log $\calC$-concave,
\item $G$ is almost surely log-concave.
\end{enumerate}
\end{thm}

The notion of a convex set can be extended as
\begin{defn}
Any measurable subset $A$ of $W$ will be called $H$-convex  if its indicator
function $1_A$ is log $H$-concave.
\end{defn}
\remark
Evidently any measurable  convex subset of $W$ is
$H$-convex. Moreover, if $A=A'$ almost surely and if $A$ is
$H$-convex, then $A'$ is also $H$-convex. 

\remark 
If $\phi$ is an $H$-convex Wiener functional, then 
the set 
$$
\{w\in W: \phi(w)\leq t\}
$$ 
is $H$-convex for any $t\in
\reals$.

We have the following result about the characterization of the
$H$-convex sets:
\begin{thm}
Assume that $A$ is an $H$-convex set, then there exists a convex set
$A'$, which is Borel measurable such that $A=A'$ almost surely.
\end{thm}
\proof 
Since, by definition, $1_A$ is a log $H$-concave Wiener functional,
from Theorem \ref{version-thm}, there exists a log-concave Wiener
functional $f_A$ such that $f_A=1_A$ almost surely. It suffices to
define $A'$ as the set  
$$
A'=\{w\in W:\,f_A(w)\geq 1\}\,.
$$
\qed
\paragraph{Example:}
Assume that $A$ is an $H$-convex subset of $W$ of positive
measure. Define $p_A$ as
$$
p_A(w)=\inf\left(|h|_H:h\in (A-w)\cap H\right)\,.
$$
Then $p_A$ is $H$-convex, hence almost surely convex (and $H$-Lipschitz 
c.f. \cite{USZ99}). Moreover, the  $\{w:\,p_A(w)\leq \alpha\}$ is an
$H$-convex set for any $\alpha\in \reals_+$. 

\section{Extensions and some applications}

\begin{defn}
Let $(e_i,i\in \NN)$ be any  complete orthonormal  basis of  $H$. We shall
denote, as before, by $w_n=\sum_{i=1}^n\delta e_i(w)\,e_i$ and
$w_n^\perp=w-w_n$.
Assume now that $F:W\to \reals\cup\{\infty\}$ is a measurable mapping with 
$\mu\{F<\infty\}>0$. 
\begin{enumerate}
\item We say that it is $a$-convex ($a\in \reals$), if the partial map
$$
w_n\to \frac{a}{2}|w_n|^2+F(w_n^\perp+w_n)
$$
is almost surely convex for any $n\geq 1$, where $|w_n|$ is the
Euclidean norm of $w_n$. 
\item We call $G$ $a$-log-concave if 
$$
w_n\to \exp\left\{-\frac{a}{2}|w_n|^2\right\}G(w_n^\perp+w_n)
$$
\end{enumerate}
is almost surely log-concave for any $n\in \NN$.
\end{defn}
\remark
$G$ is $a$-log-concave if and only if $-\log G$ is $a$-convex.
\newline

\noindent
The following theorem gives  a practical method to verify
$a$-convexity or log-concavity:
\begin{thm}
\label{a-con-thm}
Let $F:W\to \reals\cup\{\infty\}$ be a measurable map such that
$\mu\{F<\infty\}>0$. Define the map $F_a$ on $H\times W$ as
$$
F_a(h,w+h)=\frac{a}{2}|h|_H^2+F(w+h)\,.
$$
Then $F$ is $a$-convex if and only if, for any $h,k\,\in H$ and  
$\alpha,\beta\in[0,1]$ with $\alpha+\beta=1$, one has 
\begin{equation}
\label{H-cond}
F_a(\alpha h+\beta k,w+\alpha h+\beta k)\leq \alpha\, F_a(h,w+h)+\beta\,
F_a(k,w+k)
\end{equation}
$\mu$-almost surely, where the negligeable set on which the inequality 
(\ref{H-cond}) fails may depend on the choice of $h,k$ and of
$\alpha$. 

Similarly a measurable mapping $G:W\to \reals_+$ is
$a$-log-concave if and only if the map defined by 
$$
G_a(h,w+h)=\exp\left\{-\frac{a}{2}|h|_H^2\right\}G(w+h)
$$
satisfies the inequality
\begin{equation}
\label{a-cond}
G_a(\alpha h+\beta k,w+\alpha h+\beta k)\geq G_a( h,w+h)^\alpha 
G_a( k,w+k)^\beta\,,
\end{equation}
$\mu$-almost surely, where the negligeable set on which the inequality 
(\ref{a-cond}) fails may depend on the choice of $h,k$ and of
$\alpha$. 
\end{thm}
\proof 
Let us denote by $h_n$ its projection on
the vector space spanned by $\{e_1,\ldots,e_n\}$, i.e.  $h_n=\sum_{i\leq
  n}(h,e_i)_H e_i$. Then, from Theorem \ref{equiv-prop}, $F$ is
$a$-convex if and only if the map
$$
h_n\rightarrow
\frac{a}{2}\left[|w_n|^2+2(w_n,h_n)+|h_n|^2\right]+F(w+h_n)
$$ 
satisfies a convexity inequality like (\ref{H-cond}). Besides the
term $|w_n|^2$ being kept constant in this operation,   it can be
removed from the both sides of the inequality. Similarly, since
$h_n\to (w_n,h_n)$ is being affine, it also cancels from the both
sides of this inequality. Hence $a$-convexity is equivalent to
$$ 
F_a(\alpha h_n+\beta k_n,w+\alpha h_n+\beta k_n)\leq
 \alpha\, F_a(h_n,w+h_n)+\beta\,F_a(k_n,w+k_n)
$$
where $k_n$ is defined as $h_n$ from a $k\in H$. 

The second part of the theorem is obvious since $G$ is $a$-log-concave 
if and only if $-\log G$ is $a$-convex.
\qed
\begin{cor}
\label{car-cor}
\begin{enumerate}
\item Let $\hat{L}^0(\mu)$ be the space of the $\mu$-equivalence
  classes of   $\R\cup\{\infty\}$-valued
  random variables regarded as a topological semi-group under addition 
  and convergence in probability. Then $F\in \hat{L}^0(\mu)$ is
  $\beta$-convex if and only if the mapping 
$$
h\to \frac{\beta}{2}|h|_H^2+F(w+h)
$$
is a convex and continuous mapping from $H$ into $\hat{L}^0(\mu)$.
\item $F\in L^p(\mu),\,p>1$ is $\beta$-convex if and only if 
$$
E\left[\left((\beta I_H+\nabla^2F)h,h\right)_H\,\phi\right]\geq 0
$$
for any $\phi\in \DD$ positive and $h\in H$, where $\nabla^2F$ is to
be understood in the sense of the distributions $\DD'$.
\end{enumerate}
\end{cor}

\paragraph{Example:}
Note for  instance that $\sin \delta h$ with $|h|_H=1$, is a 1-convex
and that $\exp(\sin \delta h)$ is $1$-log-concave.

The following result is a direct consequence of Prekopa's theorem:
\begin{prop}
\label{c-exp-prop}
Let $G$ be an $a$-log concave Wiener functional, $a\in [0,1]$,  and
assume that $V$ is  
any sigma algebra generated by the elements of the first Wiener
chaos. Then $E[G|V]$ is again $a$-log-concave.
\end{prop}
\proof
From Corollary \ref{car-cor}, it suffices to prove the case $V$ is
generated by $\{\delta e_1,\ldots,\delta e_k\}$, where $(e_n,n\in
\NN)$ is an orthonormal basis of $H$. Let 
\beaa
w_k&=&\sum_{i\leq k}\delta e_i(w)e_i\\
z_k&=&w-w_k\\ 
z_{k,n}&=&\sum_{i=k+1}^{k+n}\delta e_i(w)e_i
\eeaa
 and let  $z_{k,n}^\perp=z_k-z_{k,n}$. Then we have 
\beaa
E[G|V]&=&\int G(z_k+w_k)d\mu(z_k)\\
&=&\lim_n\frac{1}{(2\pi)^{n/2}}\int_{\R^n}G(z_{k,n}^\perp+z_{k,n}+w_k)
e^{-\frac{|z_{k,n}|^2}{2}}dz_{k,n}\,.
\eeaa
Since 
$$
(z_{k,n},w_k)\to
\exp\left\{-\frac{1}{2}(a|w_k|^2+|z_{n,k}|^2)\right\}
G(z_{k,n}^\perp+z_{k,n}+w_k) 
$$
is almost surely log-concave, the proof follows from Prekopa's
theorem (cf.\cite{PRE71}).
\qed

\noindent
The following theorem extends Theorem  \ref{Gamma-thm} :
\begin{thm}
\label{general-O-U}
Let $G$ be an $a$-log-concave Wiener functional, where $a\in
[0,1)$. Then   $\Gamma(A)G$ is
$a$-log-concave, where $A\in L(H,H)$ (i.e. the space of bounded linear 
operators on $H$) with $\|A\|\leq 1$. In particular $P_tG$ is
$a$-log-concave for any $t\geq 0$, where $(P_t,t\geq 0)$ denotes the
Ornstein-Uhlenbeck semi-group on $W$.
\end{thm}
\proof
Let $(e_i,i\in \NN)$ be a complete, orthonormal basis of $H$, denote
by $\pi_n$ the orthogonal projection from $H$ onto the linear  space spanned 
by $\{e_1,\ldots,e_n\}$ and by $V_n$ the sigma algebra generated by
$\{\delta e_1,\ldots,\delta e_n\}$. From Proposition \ref{c-exp-prop}
and from the fact that $\Gamma(\pi_nA\pi_n)\to \Gamma(A)$ in the
strong operator topology as $n$ tends to infinity, it suffices to
prove the theorem when $W=\reals^n$. We may then assume that $G$ is
bounded and of compact support. Define $F$ as 
\beaa
G(x)&=&F(x)e^{\frac{a}{2}|x|^2}\\
    &=&F(x)\int_{\reals^n}e^{\sqrt{a}(x,\xi)}d\mu(\xi)\,.
\eeaa
From the hypothesis, $F$ is almost surely log-concave. Then, using the
notations explained in Section 2:
\beaa
\lefteqn{e^{-a\frac{|x|^2}{2}}\Gamma(A)G(x)}\\
&=&\int\int
F(A^*x+Sy)\exp\left\{-a\frac{|x|^2}{2}+\sqrt{a}(A^*x+Sy,\xi)\right\} 
d\mu(y)d\mu(\xi)\\
&=&(2\pi)^{-n}\int\int F(A^*x+Sy)\exp-\frac{\Theta(x,y,\xi)}{2}\, dyd\xi\,,
\eeaa
where
\beaa
\Theta(x,y,\xi)&=&a|x|^2-2\sqrt{a}(A^*x+Sy,\xi)+|y|^2+|\xi|^2\\
  &=&|\sqrt{a}x-A\xi|^2+|\sqrt{a}y-S\xi|^2+(1-a)|y|^2\,,
\eeaa
which is a convex function of $(x,y,\xi)$.  Hence the proof follows
from Pr\'ekopa's theorem (cf.\cite{PRE71}).
\qed

The following proposition  extends a well-known  finite dimensional
inequality (cf.\cite{Hu}):
\begin{prop}
Assume that  $f$ and $g$ are $H$-convex Wiener functionals such that 
$f\in L^p(\mu)$ and $g\in L^q(\mu)$ with  $p>1,\,p^{-1}=1-q^{-1}$. Then 
\begin{equation}
\label{cov-ineq}
E[f\,g]\geq E[f] E[g]+\left(E[\nabla f],E[\nabla g]\,\right)_H\,.
\end{equation}
\end{prop}
\proof
Define the smooth and convex functions $f_n$ and $g_n$ on $W$ by 
\beaa
P_{1/n}f&=&f_n\\
P_{1/n}g&=&g_n\,.
\eeaa
Using the fact that $P_t=e^{-tL}$, where $L$ is the number operator
$L=\delta\circ \nabla$ and the commutation relation $\nabla
P_t=e^{-t}P_t\nabla$, for any $0\leq t\leq T$,  we have 
\begin{eqnarray}
\label{yosida}
E\left[P_{T-t}f_n\,g_n\right]&=&E[P_Tf_n\,g_n]+\int_0^tE\left[LP_{T-s}f_n\,g_n\right]ds\nonumber\\
  &=&E[P_Tf_n\,g_n]+\int_0^te^{-(T-s)}E\left[\left(P_{T-s}\nabla f_n,\nabla
    g_n\right)_H\right]ds\nonumber\\ 
&=&E[P_Tf_n\,g_n]+\int_0^te^{-(T-s)}E\left[\left(P_T\nabla f_n,\nabla
g_n\right)_H\right]ds\nonumber\\
&&+e^{-2T}\int_0^t\int_0^se^{s+\tau}E\left[\left(P_{T-\tau}\nabla^2f_n,\nabla^2g_n\right)_2\right]
d\tau ds\nonumber\\
&\geq &E[P_Tf_n\,g_n]+E\left[\left(P_T\nabla f_n,\nabla g_n\right)_H\right]\,e^{-T}(e^t-1)
\end{eqnarray}
where $(\cdot,\cdot)_2$ denotes the
Hilbert-Schmidt scalar product and the inequality (\ref{yosida})
follows from the convexity of $f_n$ and $g_n$. In fact their convexity 
implies that $P_t \nabla^2 f_n$ and $\nabla^2 g_n$ are positive
operators, hence their Hilbert-Schmidt tensor product is positive.
Letting $T=t$ in the above inequality we have 
\begin{equation}
\label{yosida-2}
E[f_n\,g_n]\geq
E\left[P_Tf_n\,g_n\right]+(1-e^{-T})E\left[\left(P_T\nabla f_n,\nabla
    g_n\right)_H\right]\,. 
\end{equation}
Letting  $T\rightarrow \infty$ in (\ref{yosida-2}), we obtain, by the
ergodicity of $(P_t,t\geq 0)$, the
claimed inequality for $f_n$ and $g_n$. It suffices then to  take the
limit of this inequality  as $n$ tends to infinity.
\qed


\begin{prop}
\label{con-con}
Let $G$ be a (positive) $\gamma$-log-concave Wiener functional with $\gamma\in
[0,1]$. Then the map $h\to E[G(w+h)]$ is a log-concave mapping on $H$.
In particular, if $G$ is symmetric, i.e., if
$G(w)=G(-w)$, then 
$$
E[G(w+h)]\leq E[G]\,.
$$
\end{prop}
\proof
Without loss of generality, we may suppose that $G$ is bounded. Using
the usual notations, we have, for any $h$ in any finite dimensional
subspace $L$ of $H$, 
$$
E[G(w+h)]=\lim_n\frac{1}{(2\pi)^{n/2}}\int_{W_n}G(w_n^\perp+w_n+h)
\exp\left\{-\frac{|w_n|^2}{2}\right\}dw_n\,,
$$
from the hypothesis, the integrand is almost surely log-concave on
$W_n\times L$, from Prekopa's theorem, the integral is log-concave on
$L$, hence the limit is  also log-concave. Since $L$ is arbitrary,
the first part of the  proof follows. To prove the second part, let
$g(h)=E[G(w+h)]$, then, from the log-concavity of $g$ and symmetry of
$G$, we have 
\beaa
E[G]&=&g(0)\\
&=&g\left(1/2(h)+1/2(-h)\right)\\
          &\geq&g(h)^{1/2}g(-h)^{1/2}\\
          &=&g(h)\\
   &=&E[G(w+h)]\,.
\eeaa
\qed

\remark In fact, with a little bit more attention, we can see that
the map $h\to \exp\{\frac{1}{2}(1-\gamma)|h|^2_H\}E[G(w+h)]$ is
log-concave on $H$. 
\newline

\noindent
We have the following immediate corollary:
\begin{cor}
Assume that $A\subset W$ is an $H$-convex  and  symmetric set. Then we
have 
$$
\mu(A+h)\leq \mu(A)\,,
$$
for any $h\in H$.
\end{cor}
\proof
Since $1_A$ is log $H$-concave, the proof  follows from Proposition
\ref{con-con}.
\qed

\begin{prop}
Let $F\in L^p(\mu)$ be a  positive  log $H$-convex function. Then for
any $u\in \DD_{q,2}(H)$, we have 
$$
E_F\left[\left(\delta u-E_F[\delta u]\right)^2\right]\geq
 E_F\left[|u|_H^2+2\delta(\nabla_u u)+\trace (\nabla u\cdot \nabla u)\right]\,,
$$
where $E_F$ denotes the mathematical expectation with respect to the
probability defined as 
$$
\frac{F}{E[F]}d\mu\,.
$$
\end{prop} 
\proof Let $F_\tau$ be $P_\tau F$, where $(P_\tau,\tau\in \reals_+)$
denotes the Ornstein-Uhlenbeck semi-group. $F_\tau$ has a
modification, denoted again by the same letter, such that the mapping
$h\mapsto F_\tau(w+h)$ is real-analytic on $H$ for all $w\in W$
(cf. \cite{USZ99}). Suppose first also that $\|\nabla u\|_{_2}\in
L^\infty (\mu, H\otimes H)$ where $\|\cdot\|_{_2}$ denotes the
Hilbert-Schmidt norm. Then, for any $r>1$, there exists some $t_r>0$
such that, for any $0\leq t<t_r$,  the image of the Wiener measure
under $w\mapsto w+tu(w)$ is equivalent to $\mu$ with the Radon-Nikodym 
density $L_t\in L^r(\mu)$. 
Hence $w\mapsto F_\tau(w+tu(w))$ is a well-defined
mapping on $W$ and it is in some $L^r(\mu)$ for small $t>0$
(cf. \cite{USZ99}, Chapter~3 and  Lemma B.8.8). Besides
$t\mapsto F(w+tu(w))$ is log convex on $\reals$ 
since $F_\tau $ is log $H$-convex. Consequently $t\mapsto
E[F_\tau(w+tu(w))]$ is log convex and strictly positive. Then the
second derivative of its logarithm  at $t=0$ should be positive. This
implies immediately the claimed inequality for $\nabla u$ bounded. We
then pass to the limit with respect to $u$ in $\DD_{q,2}(H)$ and then
let $\tau \to 0$ to complete the proof.
\qed

\section{Poincar\'e and logarithmic Sobolev inequalities}
The following theorem  extends the Poincar\'e- Brascamp-Lieb 
inequality (cf.\cite{UST95}):
\begin{thm}
\label{B-L-thm}
Assume that $F$ is a Wiener functional in
$\cup_{p>1}\DD_{p,2}$ with $e^{-F}\in L^1(\mu)$ and assume also  that  there
exists a constant $\eps>0$ such that 
\begin{equation}
\label{growth-cond}
\left((I_H+\nabla^2F)h,h\right)_H\geq \eps |h|_H^2
\end{equation}
almost surely, for any $h\in H$, i.e. $F$ is $(1-\eps)$-convex.
Let us denote by $\nu_F$ the
probability measure on $(W,\calB(W))$ defined by
$$
d\nu_F=\exp\left\{-F-\log E\left[e^{-F}\right]\right\}d\mu\,.
$$ 
Then for any smooth cylindrical Wiener 
functional $\phi$, we have 
\begin{equation}
\label{poincare-ineq}
\int_W|\phi-E_{\nu_F}[\phi]|^2d\nu_F\leq \int_W\left((I_H+\nabla^2F)^{-1}\nabla
\phi,\nabla \phi\right)_Hd\nu_F\,.
\end{equation}
In particular, if $F$ is an $H$-convex Wiener functional, then the
condition (\ref{growth-cond}) is satisfied with $\eps=1$.
\end{thm}
\proof
Assume first that $W=\R^n$ and that $F$ is a smooth  function
on $\R^n$ satisfying the inequality (\ref{growth-cond}) in this
setting. Assume also for the typographical facility that
$E[e^{-F}]=1$. For any smooth function function $\phi$ on $\R^n$, we have  
\begin{equation}
\label{fin-dim}
\int_{\R^n}\left|\phi-E_{\nu_F}[\phi]\right|^2d\nu_F=\frac{1}{(2\pi)^{n/2}}
\int_{\R^n} e^{-F(x)-|x|^2/2}\left|\phi(x)-E_F[\phi]\right|^2dx\,. 
\end{equation}
The function $G(x)=F(x)+\frac{1}{2}|x|^2$ is a strictly convex smooth
function. Hence Brascamp-Lieb inequality (cf.\cite{BLI})
implies that: 
\beaa
\int_{\R^n}\left|\phi-E_{\nu_F}[\phi]\right|^2d\nu_F&\leq&\int_{\R^n}\left
(\left({\mbox{\rm 
        Hess}}\,G(x)\right)^{-1}\nabla\phi(x),\nabla\phi(x)\right)_{\R^n}d\nu_F(x)\\
&=&\int_{\R^n}\left((I_{\R^n}+\nabla^2F)^{-1}\nabla\phi,\nabla\phi\right)_
{\R^n}d\nu_F\,.
\eeaa
To prove the general case we proceed by approximation as before:
indeed let $(e_i,i\in \NN)$ be a complete, 
orthonormal basis of $H$, denote by $V_n$ the sigma algebra generated
by $\{\delta e_1,\ldots,\delta e_n\}$. Define $F_n$ as to be
$E[P_{1/n}F|V_n]$, where $P_{1/n}$ is the Ornstein-Uhlenbeck semigroup 
at $t=1/n$. Then from the martingale convergence theorem and
the fact that $V_n$ is a smooth sigma algebra, the sequence $(F_n,n\in 
\NN)$ converges to $F$ in some $\DD_{p,2}$. Moreover $F_n$ 
satisfies the hypothesis (with  a better constant in the inequality
 (\ref{growth-cond})) since $\nabla^2F_n=e^{-2/n}E[Q_n^{\otimes
   2}\nabla^2F|V_n]$, where $Q_n$ denotes the orthogonal projection
 onto the vector space spanned by $\{e_1,\ldots,e_n\}$. Besides $F_n$
 can be represented as $F_n=\theta(\delta e_1,\ldots,\delta e_n)$,
 where $\theta$ is a smooth  function on $\R^n$ satisfying 
$$
((I_{\reals^n}+\nabla^2\theta(x))y,y)_{\reals^n}\geq \eps
|y|_{\reals^n}^2\,,
$$
for any $x,y\in \reals^n$. Let $w_n={\tilde{Q}}_n(w)=\sum_{i\leq
  n}(\delta e_i)e_i$, $W_n={\tilde{P}}_n(W)$ and
$W_n^\perp=(I_W-{\tilde{Q}}_n)(W)$ as before.
Let us denote by $\nu_n$ the probability measure corresponding to
$F_n$. Let us also denote by $V_n^\perp$ the sigma algebra generated
by $\{\delta e_k,k>n\}$.  Using the finite dimensional result that we
have derived,  the Fubini theorem and the inequality $2|ab|\leq \kappa 
a^2+\frac{1}{\kappa} b^2$, for any $\kappa>0$,  we obtain 
\bea
\label{uzun-eqn}
\lefteqn{E_{\nu_{n}}\left[\left|\phi-E_{\nu_{n}}[\phi]\right|^2\right]}\nonumber\\
&=&\int_{W_n\times
  W_n^\perp}e^{-F'_n(w_n)}|\phi(w_n+w_n^\perp)-E_{\nu_n}[\phi]|^2d\mu_n(w_n)d\mu_n^\perp(w_n^\perp)\nonumber\\  
&\leq&
(1+\kappa)\int_We^{-F'_n}|\phi-E[e^{-F'_n}\phi|V_n^\perp]|^2d\mu\nonumber\\
&&+\left(1+\frac{1}{\kappa}\right)\int_We^{-F'_n}|E[e^{-F'_n}\phi|V_n^\perp]-E_{\nu_n}[\phi]|^2d\mu\nonumber\\
&\leq&
(1+\kappa)E_{\nu_n}\left[\left((I_H+\nabla^2F_n)^{-1}\nabla\phi,\nabla\phi\right)_H\right]\nonumber\\
&&+\left(1+\frac{1}{\kappa}\right)\int_We^{-F'_n}|E[e^{-F'_n}\phi|V_n^\perp]-E_{\nu_n}[\phi]|^2d\mu\,, 
\eea
where $F_n'$ denotes $F_n-\log E[e^{-F_n}]$. Since $V_n$ and
$V_n^\perp$ are independent sigma algebras, we have 
\beaa
|E[e^{-F'_n}\phi|V_n^\perp]|&=&\frac{1}{E[e^{-F_n}]}
|E[e^{-F'_n}\phi|V_n^\perp]|\\
&\leq&\frac{1}{E[e^{-F_n}]}E[e^{-F_n}|V_n^\perp]\|\phi\|_\infty\\
&=&\|\phi\|_\infty\,,
\eeaa
hence, using the triangle inequality  and the dominated convergence
theorem, we realize that the last term in (\ref{uzun-eqn}) converges
to zero as $n$ tends to infinity. Since the sequence of  operator valued random
variables $((I_H+\nabla^2F_n)^{-1},n\in \NN)$ is  essentially bounded
in the strong operator norm, we can pass to the limit on both sides
and this gives the claimed inequality with a factor $1+\kappa$, since
$\kappa>0$ is arbitrary, the  proof is completed.
\qed


\remark  
Let $T:W\to W$ be a shift  defined as $T(w)=w+u(w)$, where $u:W\to H$
is a measurable map satisfying $(u(w+h)-u(w),h)_H\geq -\eps|h|^2$. In 
\cite{USZ96} and in \cite{USZ99}, Chapter 6, we have studied such
transformations, called $\eps$-monotone shifts. Here  the hypothesis of Theorem
\ref{B-L-thm} says that the shift $T=I_W+\nabla F$ is $\eps$-monotone.

The Sobolev regularity hypothesis can be omitted if we are after a
Poincar\'e inequality with another  constant:
\begin{thm}
Assume that $F\in \cup_{p>1}L^p(\mu)$ with
$E\left[e^{-F}\right]$ is finite and that, for some constant $\eps>0$, 
$$
E\left[\left((I_H+\nabla^2F)h,h\right)_H\,\psi\right]\geq
\eps\,|h|_H^2E[\psi]\,,
$$
for any $h\in H$ and positive test function $\psi\in\DD$, where
$\nabla^2F$ denotes the second order derivative in the sense of the
distributions. Then we have  
\begin{equation}
\label{poincare-2-ineq}
E_{\nu_F}\left[|\phi-E_F[\phi]|^2\right]\leq
\frac{1}{\eps}E_{\nu_F}[|\nabla\phi|_H^2]
\end{equation}
for any cylindrical Wiener functional $\phi$. In particular, if $F$ is 
$H$-convex, then we can take $\eps=1$.
\end{thm}
\proof
Let $F_t$ be defined as $P_tF$, where $P_t$ denotes the
Ornstein-Uhlenbeck semigroup. Then $F_t$ satisfies the hypothesis of
Theorem \ref{B-L-thm}, hence we have 
$$
E_{\nu_{F_t}}\left[\left|\phi-E_{F_t}[\phi]\right|^2\right]\leq
\frac{1}{\eps}E_{\nu_{F_t}}\left[|\nabla\phi|_H^2\right]
$$
for any $t>0$. The claim follows when we take the limits of both sides 
as $t\to 0$.
\qed
\paragraph{\bf Example:}
Let $F(w)=\|w\|+\frac{1}{2}\sin(\delta h)$ with $|h|_H\leq 1$, where
$\|\cdot\|$ denotes the norm of the Banach 
space $W$. Then in general $F$ is not in $\cup_{p>1}\DD_{p,2}$,
however the Poincar\'e  inequality (\ref{poincare-2-ineq}) holds with
$\eps=1/2$. 
\begin{thm}
\label{log-sob-thm}
Assume that $F$ is a Wiener functional in
$\cup_{p>1}\DD_{p,2}$ with $E[\exp-F]<\infty$. Assume  that  there
exists a constant $\eps>0$ such that 
\begin{equation}
\label{growth-cond-1}
\left((I_H+\nabla^2F)h,h\right)_H\geq \eps |h|_H^2
\end{equation}
almost surely, for any $h\in H$.
Let us denote by $\nu_F$ the
probability measure on $(W,\calB(W))$ defined by
$$
d\nu_F=\exp\left\{-F-\log E\left[e^{-F}\right]\right\}d\mu\,.
$$ 
Then for any smooth cylindrical Wiener 
functional $\phi$, we have 
\begin{equation}
\label{log-sob-ineq}
E_{\nu_F}\left[\phi^2\left\{\log \phi^2-\log \|\phi\|_{L^2(\nu_F)}^2\right\}\right]\leq
\frac{2}{\eps} E_{\nu_F}\left[|\nabla \phi|_H^2\right]\,.
\end{equation}
In particular, if $F$ is an $H$-convex Wiener functional, then the
condition (\ref{growth-cond-1}) is satisfied with $\eps=1$.
\end{thm}
\proof 
We shall proceed as in the proof of Theorem \ref{B-L-thm}. Assume then 
that $W=\R^n$ and that $F$ is a smooth function satisfying the
inequality (\ref{growth-cond-1})  in this frame. In this case it is
immediate to see that  function $G(x)=\frac{1}{2}|x|^2+F(x)$ satisfies 
the Bakry-Emery condition (cf.\cite{B-E}, \cite{D-S}), which is known
as a sufficient condition for the inequality (\ref{log-sob-ineq}). For 
the infinite dimensional case we define as in the proof of Theorem
\ref{B-L-thm}, $F_n,\nu_n,V_n,V_n^\perp$. Then, denoting by $E_n$ the
expectation with respect to the probability $\exp\{-F'_n\}d\mu$, where
$F_n'=F_n-\log E[e^{-F_n}]$, we have  
\bea
\label{l-s-1}
\lefteqn{E_n\left[\phi^2\left\{\log \phi^2-\log
    \|\phi\|_{L^2(\nu_F)}^2\right\}\right]}\nonumber\\
&=&E_n\left[\phi^2\left\{\log \phi^2-\log
    E[e^{-F_n'}\phi^2|V_n^\perp]\right\}\right]\nonumber\\
&&+E_n\left[\phi^2\left\{\log E[e^{-F_n'}\phi^2|V_n^\perp]-\log E_n[\phi^2]\right\}\right]\nonumber\\
&\leq&\frac{2}{\eps}E_n\left[|\nabla
  \phi|_H^2\right]+E_n\left[\phi^2\left\{\log
    E[e^{-F_n'}\phi^2|V_n^\perp]-\log E_n[\phi^2]\right\}\right]\,, 
\eea
where we have used, as in the proof of Theorem \ref{B-L-thm}, the
finite dimensional log-Sobolev inequality to obtain the inequality  
(\ref{l-s-1}). 
Since in the above inequalities everything is squared, we can assume
that $\phi$ is positive, and adding a constant $\kappa>0$, we can also 
replace $\phi $ with  $\phi_\kappa=\phi+ \kappa$. Again by the
independance of $V_n$ and $V_n^\perp$, we  can pass to the limit with
respect to $n$ in
the inequality (\ref{l-s-1}) for $\phi=\phi_\kappa$ to obtain 
$$
E_{\nu_F}\left[\phi_\kappa^2\left\{\log \phi_\kappa^2-\log \|\phi_\kappa\|_{L^2(\nu_F)}^2\right\}\right]\leq
\frac{2}{\eps} E_{\nu_F}\left[|\nabla \phi_\kappa|_H^2\right]\,.
$$
To complete the proof it suffices to pass to the limit as $\kappa\to
0$.
\qed

The following theorem fully   extends Theorem \ref{log-sob-thm} and it 
is useful for the applications:
\begin{thm}
\label{l-s-3}
Assume that $G$ is a (positive) $\gamma$-log-concave Wiener functional for
some $\gamma\in [0,1)$ with $E[G]<\infty$. Let us denote by $E_G[\cdot\,]$
the expectation with respect to the probability measure defined by 
$$
d\nu_G=\frac{G}{E[G]}d\mu\,.
$$
Then we have 
\begin{equation}
\label{l-s-2}
E_G\left[\phi^2\left\{\log \phi^2-\log E_G[\phi^2]\right\}\right]\leq
\frac{2}{1-\gamma}E_G[|\nabla \phi|_H^2]\,,
\end{equation}
for any cylindrical Wiener functional $\phi$.
\end{thm}
\proof
Since $G\wedge c$, $c>0$,  is again $\gamma$-log-concave, we may
suppose without loss of generality that $G$ is bounded. Let now
$(e_i,i\in \NN)$ be a 
complete, orthonormal basis for $H$, denote by $V_n$ the sigma algebra 
generated by $\{\delta e_1,\ldots,\delta e_n\}$. Define $G_n$ as to be 
$E[P_{1/n}G|V_n]$. From Proposition \ref{c-exp-prop} and Theorem
\ref{general-O-U}, $G_n$ is again a  $\gamma$-log-concave, strictly positive
Wiener functional. It can be represented  as 
$$
G_n(w)=g_n(\delta e_1,\ldots,\delta e_n)
$$
and due to the Sobolev embedding  theorem, after a modification on a
set of zero Lebesgue measure, we can assume that  $g_n$ is  a smooth
function on $\reals^n$. Since it is 
strictly positive, it is of the form $e^{-f_n}$, where $f_n$ is a
smooth, $\gamma$-convex function. It follows then from Theorem
\ref{log-sob-thm} that the inequality (\ref{l-s-2}) holds when we
replace $G$ by $G_n$, then the proof follows by taking the limits  of
both  sides as $n\to \infty$.
\qed
\paragraph{Example:}
Assume that $A$ is a measurable subset of $W$ and let $H$ be a
measurable Wiener functional with values in
$\reals\cup\{\infty\}$. If $G$ defined by  $G=1_A\,H$ is
$\gamma$-log-concave with $\gamma\in [0,1)$, then the hypothesis of Theorem
\ref{l-s-3} are satisfied.

\begin{defn}
Let $T\in \DD'$ be a positive distribution. We say that it is
$a$-log-concave if $P_tT$ is an $a$-log-concave Wiener functional. If
$a=0$, then we call $T$ simply log-concave.
\end{defn}
\remark
It is well-known that (cf. for example \cite{UST95}), to any positive
distribution on $W$, it 
corresponds a positive Radon measure $\nu_T$ such that 
$$
<T,\phi>=\int_W {\tilde{\phi}}(w)d\nu_T(w)
$$
 for any $\phi\in \DD$, where $\tilde{\phi}$ represents a
 quasi-continuous version of $\phi$. 
\paragraph{Example:}
Let $(w_t,t\in [0,1])$ be the one-dimensional Wiener process and
denote by  $p_\tau$ the heat kernel on $\reals$.  Then the 
distribution defined as $\veps_0(w_1)=\lim_{\tau\to 0}p_\tau(w_1)$ is
log-concave, where $\veps_0$ denotes the Dirac measure at zero.

The following result is a Corollary of Theorem \ref{l-s-3}:

\begin{thm}
\label{dist-case}
Assume that $T\in \DD'$ is a positive, $\beta$-log-concave distribution
with $\beta\in [0,1)$. Let $\gamma$ be the probability  Radon  measure
defined by
$$
\gamma=\frac{\nu_T}{<T,1>}\,.
$$
Then we have 
\begin{equation}
\label{l-s-4}
E_\gamma\left[\phi^2\left\{\log \phi^2-\log E_\gamma[\phi^2]\right\}\right]\leq
\frac{2}{1-\beta}E_\gamma[|\nabla \phi|_H^2]\,,
\end{equation} 
for any smooth cylindrical function $\phi:W\to \reals$.
\end{thm}

Here is  an application of this result:
\begin{prop}
\label{non-deg}
Let $F$ be a Wiener functional in $\DD_{r,2}$ for some $r>1$. Suppose
that it is $p$-non-degenerate in the sense that 
\begin{equation}
\label{non-deg-con}
\delta\left\{\frac{|\nabla F|_H^2}{|F|^2}\phi\right\}\in L^p(\mu)
\end{equation}
for any $\phi\in \DD$, for some $p>1$. Assume furthermore that, for
some $x_0\in \reals$ and $a\in [0,1)$, 
\begin{equation}
\label{x0-cond}
(F-x_0)\nabla^2F+\nabla F\otimes \nabla F\geq -a I_H
\end{equation}
almost surely. Then we have 
$$
E\left[\phi^2\left\{\log \phi^2-\log
    E\left[\phi^2|F=x_0\right]\right\}|F=x_0\right]\leq
\frac{2}{1-a}E\left[|\nabla \phi|_H^2|F=x_0\right]
$$
for any smooth cylindrical $\phi$.
\end{prop}
\proof
Note that the non-degeneracy  hypothesis (\ref{non-deg-con}) implies
the existence of a continuous 
density of the law of $F$ with respect to the Lebesgue
measure (cf. \cite{MAL97} and the references there). Moreover it
implies also the fact that  
$$
\lim_{\tau\to 0}p_\tau(F-x_0)=\veps_{x_0}(F)\,,
$$
in $\DD'$, where $\veps_{x_0}$ denotes the Dirac measure at $x_0$ and
$p_\tau$ is  the heat kernel on $\reals$. The
inequality (\ref{x0-cond}) implies that the distribution 
defined by
$$
\phi\to E[\phi|F=x_0]=\frac{<\veps_{x_0}(F),\phi>}{<\veps_{x_0}(F),1>}
$$
is $a$-log-concave, hence the conclusion follows from Theorem
\ref{dist-case}.
\qed
\section{Change of variables formula and log-Sobolev inequality}

In this section we shall derive a different kind of logarithmic
Sobolev inequality using the change of variables formula for the
monotone shifts studied in \cite{USZ96} and in more detail in
\cite{USZ99}. An analogous approach to derive log-Sobolev-type
inequalities  using the Girsanov theorem has 
been employed  in \cite{UST00}. 

\begin{thm}
\label{log-sob-tordu}
Suppose that $F\in L^p(\mu)$, for some $p>1$, is an $a$-convex Wiener
functional, $a\in [0,1)$ with $E[F]=0$. Assume that 
\begin{equation}
\label{growth2-cond}
E\left[\exp\left\{ c\,\|\nabla^2 L^{-1}F\|_{_2}^2\right\}\right]<\infty\,,
\end{equation}
for some 
$$
c>\frac{2+(1-a)}{2(1-a)}\,,
$$
where $\|\cdot\|_{_2}$ denotes the Hilbert-Schmidt norm on $H\otimes
H$ and $L^{-1}F=\int_{\R_+}P_tF\,dt$. Denote by $\nu$ the 
probability measure defined by 
$$
d\nu=\La\,d\mu\,,
$$
where 
$$
\La=\dett(I_H+\nabla^2L^{-1}F)\exp\left\{-F-\frac{1}{2}|\nabla
  L^{-1}F|_H^2\right\}\,
$$
and $\dett(I_H+\nabla^2L^{-1}F)$ denotes the modified Carleman-Fredholm
determinant. Then we have 
\begin{equation}
\label{eqn-tordu}
E_\nu\left[f^2\log\left(\frac{f^2}{\|f\|^2_{L^2(\nu)}}\right)\right]\leq
2 E_\nu\left[|(I_H+\nabla^2 L^{-1}F)^{-1}\nabla f|_H^2\right] 
\end{equation}
and 
\begin{equation}
\label{poincare-tordu}
E_\nu[|f-E_\nu[f]|^2]\leq E_\nu\left[|(I_H+\nabla^2
  L^{-1}F)^{-1}\nabla f|_H^2\right] 
\end{equation}
for any smooth, cylindrical $f$.
\end{thm}
\proof
Let $F_n=E[P_{1/n}F|V_n]$, where $V_n$ is the sigma algebra generated
by $\{\delta e_1,\ldots,\delta e_n\}$ and let  $(e_n,n\in \NN)$ be  a complete,
orthonormal basis of $H$. Define $\xi_n$ by $\nabla L^{-1}F_n$, then
$\xi_n$ is $(1-a)$-strongly monotone (cf. \cite{USZ99} or \cite{USZ96}) and
smooth. Consequently, the shift $T_n:W\to W$, defined by
$T_n(w)=w+\xi_n(w)$ is a bijection of $W$ (cf.\cite{USZ99}, Corollary 
6.4.1), whose inverse is of the form $S_n=I_W+\eta_n$, where
$\eta_n(w)=g_n(\delta e_1,\ldots,\delta e_n)$ such that
$g_n:\reals^n\to \reals^n$ is a smooth function. Moreover the images of 
$\mu$ under $T_n$ and $S_n$, denoted by $T_n^*\mu$ and $S_n^*\mu$
respectively,  are equivalent to $\mu$ and we have 
\beaa
\frac{dS_n^*\mu}{d\mu}&=&\La_n\\
\frac{dT_n^*\mu}{d\mu}&=&L_n
\eeaa
where 
\beaa
\La_n&=&\dett(I_H+\nabla \xi_n)\exp\left\{-\delta
  \xi_n-\frac{1}{2}|\xi_n|_H^2\right\}\\
L_n&=&\dett(I_H+\nabla \eta_n)\exp\left\{-\delta
  \eta_n-\frac{1}{2}|\eta_n|_H^2\right\}\,.
\eeaa
The hypothesis (\ref{growth2-cond}) implies the uniform integrability
of the densities $(\La_n,n\geq 1)$ and $(L_n,n\geq 1)$
(cf. \cite{USZ96,USZ99}). 
For any probability $P$ on $(W,\calB(W))$  and any positive, measurable
function $f$,  define  $\calH_P(f)$ as 
\begin{equation}
\label{entropy}
\calH_P(f)=f(\log f-\log E_P[f]).
\end{equation}
Using    the logarithmic Sobolev inequality of L. Gross  for $\mu$
(cf.\cite{LG}) and the relation
$$
(I_H+\nabla \eta_n)\circ T_n=(I_H+\nabla \xi_n)^{-1}\,,
$$
we have 
\begin{eqnarray}
\label{neweq}
E[\La_n \calH_{\La_nd\mu}(f^2)]&=&E[\calH_\mu(f^2\circ S_n)]\nonumber\\
&\leq& 2E[|\nabla (f\circ S_n)|_H^2]\nonumber\\
&=&2E[|(I_H+\nabla \eta_n)\nabla f\circ S_n|_H^2]\nonumber\\
&=&2E[\La_n|(I_H+\nabla \xi_n)^{-1}\nabla f|_H^2]\,.
\end{eqnarray}
 It follows by
the $a$-convexity  of $F$  
that 
$$
\|(I_H+\nabla \xi_n)^{-1}\|\leq \frac{1}{1-a}
$$ 
almost surely for any $n\geq 1$, where $\|\cdot\|$ denotes the operator norm. 
Since the sequence $(\La_n,n\in \NN)$ is uniformly integrable, 
 the limit of (\ref{neweq}) exists  in $L^1(\mu)$ and the proof of
 (\ref{eqn-tordu}) follows. The proof of the inequality
 (\ref{poincare-tordu}) is now trivial. 
\qed

\begin{cor}
Assume that $F$ satisfies the hypothesis of Theorem
\ref{log-sob-tordu}. Let $Z$ be the functional defined by  
$$
Z=\dett(I_H+\nabla^2 L^{-1}F)\exp\frac{1}{2}|\nabla
L^{-1}F|_H^2 
$$
and assume that $Z,\,Z^{-1}\in L^\infty(\mu)$. 
Then we have 
\begin{equation}
\label{last-ineq}
E\left[e^{-F}f^2\log\left\{\frac{f^2}{E[e^{-F}f^2]}\right\}\right]\leq
2K 
E\left[e^{-F}\left|(I_H+\nabla^2 L^{-1}F)^{-1}\nabla f\right|_H^2\right]
\end{equation}
and 
\begin{equation}
\label{last-poincare}
E\left[e^{-F}\left|f-E[e^{-F}f]\right|^2\right]\leq K E\left[e^{-F}
\left|(I_H+\nabla^2 L^{-1}F)^{-1}\nabla f\right|_H^2\right]
\end{equation}
for any smooth, cylindrical $f$, where
$K=\|Z\|_{L^\infty(\mu)}\|Z^{-1}\|_{L^\infty(\mu)}$. 
\end{cor}
\proof
Using the identity remarked by  Holley and
Stroock (cf. \cite{H-S}, p.1183)
$$
E_P\left[\calH_P(f^2)\right]=\inf_{x>0}
E_P\left[f^2\log\left(\frac{f^2}{x}\right)-(f^2-x)\right]\,, 
$$
where $P$ is an arbitrary probability measure, and $\calH$ is defined
by the relation  (\ref{entropy}),  
we see that the inequality (\ref{last-ineq}) follows from Theorem
\ref{log-sob-tordu} and  the inequality (\ref{last-poincare}) is trivial.
\qed

\remark
If $F$ is $H$-convex, then $\dett(I_H+\nabla^2L^{-1}F)\geq 1$ almost
surely. Hence in this case it suffices to assume that
$\dett(I_H+\nabla^2L^{-1}F)\in L^\infty(\mu)$ and that
$|\nabla L^{-1}F|_H\in L^\infty(\mu)$.


{\footnotesize{
\begin{tabular}{ll}
 D. Feyel & A.S. \"Ust\"unel\\
Universit\'e d'Evry-Val-d'Essone, \hspace{1.5cm }
& ENST,   D\'ept. R\'eseaux,\\
D\'ept. de Math\'ematiques, & 46, rue Barrault,  \\
91025 Evry Cedex&75013 Paris,  \\
France& France\\
 feyel@math.univ-evry.fr & ustunel@enst.fr
\end{tabular}}

\end{document}